\documentclass[11pt,a4paper]{amsart}
\usepackage{amsmath}
\usepackage{verbatim}
\usepackage{upref}
\usepackage{accents}
\usepackage{amsfonts,amssymb}
\usepackage{youngtab}
\makeatletter
\calclayout
\numberwithin{equation}{section}

\newcommand{\N}{\mathbb{N}}
\newcommand{\Z}{\mathbb{Z}}
\newcommand{\F}{F\langle X|\Z_2\rangle}
\newcommand{\f}{F\langle X\rangle}

\newcommand{\T}{\overline{t}}
\newcommand{\m}{m_{\lambda,\mu}}

\newcommand{\proofend}{\hfill{\mbox{$\Box$}}}
\newcommand{\proofbegin}{\noindent{\sc Proof.\ }}
\renewcommand{\proof}{\proofbegin}
\renewcommand{\endproof}{\proofend}

\begin{document}\title[Cocharacters of $UT_2(E)$]{Ordinary and $\Z_2$-graded Cocharacters of $UT_2(E)$}\author{Lucio Centrone}\address{Dipartimento di Matematica, Universita' degli Studi di Bari, Bari, via Orabona 4, 70125}\email{centrone@dm.uniba.it}\maketitle
\begin{abstract} Let $E$ be the infinite dimensional Grassmann algebra over a field $F$ of characteristic 0. In this paper we compute the ordinary and the $\Z_2$-graded cocharacters of the algebra of $2\times2$ upper triangular matrices with coefficients in $E$, using the tool of proper Hilbert series. \end{abstract}

\section{Introduction}
Varieties of associative algebras over a field $F$ are in bijective correspondence with the ideals of the free
associative algebra $F\langle X\rangle=F\langle x_1,x_2,\ldots\rangle$ that is invariant under all the endomorphisms of $\f$. Such
ideals are called "$T$-ideals" and they are the ideals of the polynomial identities satisfied by
any algebra of the variety. For the study of the $T$-ideals over a field of characteristic zero, a
fundamental tool is given by the representation theory of the linear and symmetric groups.
Moreover, the theorems of Kemer about the classification of the $T$-ideals of $\f$ show that
the notion of grading of an algebra defined by a group is another key ingredient for such
study. In particular, one has that any proper $T$-ideal of $\f$ is the ideal of the polynomial
identities satisfied by the Grassmann envelope of a suitable $\Z_2$-graded algebra (also called
superalgebra) of finite dimension. The work of Giambruno and Zaicev \cite{Giam02}
has contributed to clarify why the notion of PI-exponent is crucial for a classification of
the $T$-ideals in terms of growth of the sequence of their codimensions. Recall that the
$n$-th codimension of a $T$-ideal is defined as the degree of the representation of the group
$S_n$ on the vector space of the multilinear polynomials of $\f$ of degree $n$ modulo the
considered $T$-ideal. In \cite{Giam02} the authors prove that the minimal varieties with respect to
a fixed exponent are determined by the $T$-ideals of the Grassmann envelope of the so called
"minimal superalgebras". Over an algebraic closed field, such superalgebras can be
realized as graded subalgebras of block-triangular matrix algebras equipped with a suitable
$\Z_2$-grading. Precisely, the blocks along the main diagonal are simple superalgebras of finite
dimension. Then, by the Theorem of Lewin \cite{Lew} one has that the $T$-ideals of the identities
satisfied by the minimal superalgebras and their Grassmann envelopes are products of the
$T$-ideals corresponding to the diagonal blocks. Such results allow hence to solve in the
positive a conjecture due to Drensky \cite{Dr01,Dr02} about the factorability of the $T$-ideals of minimal
varieties as a product of verbally prime $T$-ideals. Moreover, Berele and Regev \cite{Ber} proved
a formula that relates the sequence of ordinary cocharacters of a product of $T$-ideals to the
sequences of cocharacters of these ideals.
Recently, in \cite{Divi01} Di Vincenzo and La Scala introduce the notion of $G$-regularity of an algebra, if $G$ is a finite abelian group, and they
prove that $T_G(R) = T_G(A)T_G(B)$,
provided that at least one of the algebras $A,B$ is G-regular. They also proved for suitable groups $G$ and for $A,B$ that the $G$-regularity of $A$ or $B$ is a necessary condition for the ideal $T_G(R)$ to be
factorable. They proved a formula allowing us to compute
the sequence of graded cocharacters of a superalgebra $R$ such that $T_2(R) = T_2(A)T_2(B)$
starting from the corresponding sequences of $A$ and $B$.

The varieties with exponent two were characterized by Giambruno and Zaicev in \cite{Giam03}. They listed five algebras such that a variety has exponent less then or equal to 2 if and only if it does not contain any of $E$ and $UT_2(F)$. These algebras generate all the possible minimal varieties of exponent strictly greater than two. One of these algebras is $G=\left(
                                              \begin{array}{cc}
                                                E & E \\
                                                0 & E^0 \\
                                              \end{array}
                                            \right),$
 where $E=E^0+E^1$ is the natural $\Z_2$-grading of $E$. In this paper, we compute firstly the Hilbert Series (and so the cocharacters) of $G$, then intends to contribute to this line of research by studying the ordinary and the $\Z_2$-graded cocharacters of  $$UT_2(E)=\left(
                                              \begin{array}{cc}
                                                E & E \\
                                                0 & E \\
                                              \end{array}
                                            \right),$$
that is a minimal algebra respect to its exponent, using the fundamental tool of proper Hilbert series.

\section{Ordinary and $\Z_2$-graded structure}
All fields in this paper are assumed to be of characteristic 0.

Let $F$ be a field and $A$ be an associative
$F$-algebra. Let now $X=\{x_1, x_2, \ldots\}$ be a countable set of variables. We denote by $\f$ the free
associative algebra generated by $X$ and by $T(A)$ the intersection of the kernels of all
 homomorphisms $\f\rightarrow A$. Then $T(A)$ is a two-sided ideal of $\f$ and its elements are called polynomial identities of the algebra $A$. If $T(A)$ is not trivial, $A$ is said to be a \emph{polynomial identity algebra} (P.I. algebra). Note that
$T(A)$ is stable under the action of any endomorphism of the algebra $\f$.
Any ideal of $\f$ which verifies such property is said to be a $T$-ideal. Clearly,
any $T$-ideal $I$ is the ideal of the polynomial identities of the algebra
$\f/I$ . Note also that for a P.I. algebra $A$, the quotient algebra $\f/T(A)$ is the
relatively free algebra for the variety of algebras generated by $A$.\vspace{3mm}

\newtheorem{definition}{Definition}[section]\begin{definition} \textrm{For $n\in\N$, the vector space $$V_n:=\text{\rm span}_F\left<x_{\sigma(1)}x_{\sigma(2)}\cdots x_{\sigma(n)}|\sigma\in S_n,\ x_i\in X\right>$$ is called the space of \emph{multilinear polynomials of degree} $n$.}\end{definition}\vspace{3mm}

Since the characteristic of the ground field $F$ is zero, a standard process of multilinearization shows that $T(A)$ is generated, as a $T$-ideal, by the subspaces $V_n\cap T(A).$ Actually, it is more efficient to study the factor space $$V_n(A)={V_n}/({V_n\cap T(A)}).$$ An effective tool to this end is provided by the representation theory of the symmetric group.
Indeed, one can notice that $V_n$ is an $S_n$-module with respect to the natural left action, and $V_n\cap T(A)$ is an $S_n$-submodule, hence the factor space $V_n(A)$ is an $S_n$-module, too. We shall denote by $\chi_n(A)$ its character, called the $n$\emph{-th cocharacter of A}. For a more detailed account about the representation theory of symmetric group, we remand to \cite{Sag}.\vspace{3mm}

The \emph{commutator} of $a,b\in A$ is the Lie product $[a,b]:=ab-ba.$ One defines inductively higher (left-normed) commutators by setting $[a_1,\cdots,a_n]:=[[a_1,\ldots,a_{n-1}],a_n],$ for any $n\geq2.$ By the Poincar$\grave{e}$-Birkhoff-Witt theorem, $\f$ has a basis $\left\{x_1^{s_1}\cdots x_r^{s_r}u_1^{m_1}\cdots u_n^{m_n}|s_i,m_j\geq0,\ r,n\in\N\right\},$ where $u_1,u_2,\ldots$ are higher commutators. We denote by $B(X)$ the unitary subalgebra of $\f$ generated by commutators, called the algebra of \emph{proper polynomials}. It is well known that $B(X)\cap T(A)$ generates the whole $T(A)$ as a $T$-ideal. Let us denote $B(A):=B(X)/B(X)\cap T(A).$ We shall denote $\Gamma_n$ the set of multilinear polynomials of $V_n$ which are proper. It is not difficult to see that $\Gamma_n$ is a left $S_n$-submodule of $V_n$ and the same holds for $\Gamma_n\cap T(A).$ Hence the factor module $$\Gamma_n(A)={\Gamma_n}/({\Gamma_n\cap T(A)})$$ is an $S_n$-submodule of $V_n(A)$. We shall denote by $\xi_n$ its character ($n$-\emph{th proper cocharacter of A}). For a more detailed account about proper cocharacters we refer to the book of Drensky (\cite{Dr03}, Chapters 4, 12).\vspace{3mm}

The following result of Drensky relates the ordinary cocharacters of a P.I. algebra $A$ with the proper cocharacters of $A$:

\newtheorem{prop}{Proposition}[section]
\begin{prop}
(Drensky, \cite{Dr03}, Theorem 12.5.4) Let A be a P.I. algebra and let $\chi_n(A)=\sum_{\lambda\vdash n}m_{\lambda}(A)\chi_{\lambda}$ its n-th cocharacter. Let $\xi_p(A)=\sum_{\nu\vdash p}k_{\nu}(A)\chi_{\nu}$ its p-th proper cocharacter, then $$m_{\lambda}(A)=\sum_{\nu\in S}k_{\nu}(A),$$ where $S=\{\nu=(\nu_1,\ldots,\nu_n)\mid\lambda_1\geq\nu_1\geq\lambda_2\geq\nu_2\geq\cdots\geq\lambda_n\geq\nu_n\}.$
\end{prop}\vspace{3mm}

We say that $A$ is a $\Z_2$-graded algebra if $A =A^0\oplus A^1$, where $A^0,A^1\subseteq A$
are subspaces and $A_gA_h\subseteq A_{g+h}$ holds for $g,h\in \Z_2$. The subspace $A_g$ is called
the homogeneous component of $A$ of degree $g$. We say that the elements $a\in A_g$ are
homogeneous of degree $g$ and we denote their degrees as $|a|=g$. One defines $\Z_2$-graded:
subspaces of $A$, $A$-modules, homomorphisms and so on, in a standard way, see for
example [1].

We denote by $\F$ the free
associative algebra generated by $X=Y\cup Z,$ where $Y=\{y_1,y_2,\ldots\}$ and $Z=\{z_1,z_2,\ldots\}$ are two countable sets of disjoint variables. Given a map $|\ |:X\rightarrow \Z_2$, we can define a $\Z_2$-grading on
$\F$ if we set $|w| = |x_{j_1}|+\cdots+|x_{j_n}|$ for any monomial $w = x_{j_1}\cdots x_{j_n}\in\F$. Then,
the homogeneous component $\F_g\subseteq\F$ is the subspace spanned by all monomials of
degree $g$. Because $\Z_2$ is a finite group, we assume that the fibers of the map $|\ |$ are all infinite.
If $A$ is a $\Z_2$-graded algebra, we denote by $T_{\Z_2}(A)$ the intersection of the kernels of all
$\Z_2$-graded homomorphisms $\F\rightarrow A$. Then $T_{\Z_2}(A)$ is a graded two-sided ideal of $\F$ and its elements are called $\Z_2$-graded polynomial identities of the algebra $A$. Note that
$T_{\Z_2}(A)$ is stable under the action of any $\Z_2$-graded endomorphism of the algebra $\F$.
Any $\Z_2$-graded ideal of $\F$ which verifies such property is said to be a $T_{\Z_2}$-ideal. Clearly,
any $T_{\Z_2}$-ideal $I$ is the ideal of the $\Z_2$-graded polynomial identities of the graded algebra
$\F/I$ . Note also that for a $\Z_2$-graded algebra $A$, the quotient algebra $\F/T_{\Z_2}(A)$ is the
relatively free algebra for the variety of graded algebras generated by $A$.\vspace{3mm}

\begin{definition} For $n\in\N$, the vector space $$V_n^{\Z_2}:=\text{\rm span}_F\left<x_{\sigma(1)}x_{\sigma(2)}\cdots x_{\sigma(n)}|\sigma\in S_n,\ x_i=y_i\text{ or }x_i=z_i\right>$$ is called the space of $\Z_2$-graded multilinear polynomials of degree $n$.\end{definition}\vspace{3mm}
Since the charcteristic of the ground field $F$ is zero, a standard process of multilinearization shows that $T_{\Z_2}(A)$ is generated, as a $T_{\Z_2}$-ideal, by the subspaces $V_n^{\Z_2}\cap T_{\Z_2}(A).$ Actually, it is more efficient to study the factor space $$V_n^{\Z_2}(A)={V_n^{\Z_2}}/({V_n^{\Z_2}\cap T_{\Z_2}(A)}).$$

One can notice again that $V_n^{\Z_2}$ is an $S_n$-module with respect to the natural left action, and $V_n^{\Z_2}\cap T_{\Z_2}(A)$ is an $S_n$-submodule, hence the factor space $V_n^{\Z_2}(A)$ is an $S_n$-module, too. We shall denote by $\chi_n^{\Z_2}(A)$ its character, called the $n$\emph{-th $\Z_2$-graded cocharacter of A}.

The study of the structure of $V_n^{\Z_2}(A)$ can be simplified by considering "smaller" spaces of multilinear polynomials. To be more precise, for fixed $l,\ m$, set $$V_{l,m}:=\textrm{span}_F\left<w\ \textrm{monomials of $V_{l+m}^{\Z_2} $}|y_1,\ldots,y_l,z_{l+1},\ldots,z_{l+m}\textrm{ occur in $w$}\right>.$$
Setting $n:=l+m$, and $S_l\times S_m=Sym(\{1,\ldots,l\})\times Sym(\{l+1,\ldots,l+m\})\leq S_n,$
the space $V_{l,m}$ is an $S_l\times S_m$-module, and the subspace $V_{l,m}\cap T_{\Z_2}(A)$ is a submodule. Therefore one can form the factor $S_l\times S_m$-module $$V_{l,m}(A):={V_{l,m}}/({V_{l,m}\cap T_{\Z_2}(A)}).$$ We shall denote by $\chi_{l,m}(A)$ its $S_l\times S_m$-character.

Now we give a small account on the representation theory of the groups $S_l\times S_m$ $(l+m=n)$. The irreducible $S_l\times S_m$-characters are in "one to one" correspondence with the pairs of partitions $(\lambda,\mu)$ of $l$ and $m$ respectively; in this case we write $\lambda\vdash l,\ \mu\vdash m,$ and $|\lambda|=l,\ |\mu|=m.$ More precisely, if $\chi_{\nu}$ denotes the irreducible $S_{|\nu|}$-character associated to the partition $\nu$, then the irreducible $S_l\times S_m$-character associated to $(\lambda,\mu)$ is $\chi_{\lambda,\mu}=\chi_{\lambda}\otimes\chi_{\mu}.$\vspace{3mm}

\section{Hilbert series and proper Hilbert series of P.I. algebras}

In this section we talk about a tool used in the study of algebras in general, the Hilbert series. In particular, we emphasize the relationship between the Hilbert series and the combinatorial method of the Littlewood-Richardson rule.\vspace{3mm}

We say that the vector space $V$ is $\Z_m$-graded if $V=\bigoplus_{(n_1,\ldots,n_m)\in\Z_m}V^{(n_1,\ldots,n_m)},$ and $V^{(n_1,\ldots,n_m)}$ is a vector subspace of $V$.
\begin{definition} Let $V=\sum_{n\in\Z}V^{(n_1,\ldots,n_m)}$ be a $\Z^m$-graded vector space and let $\dim_FV^{(n_1,\ldots,n_m)}<\infty.$ The formal power series $$H(V,t_1,\ldots,t_m)=\sum_n\dim_FV^{(n_1,\ldots,n_m)}t_1^{n_1}\cdots t_m^{n_m}$$ is called the \emph{Hilbert Series of V in the variables $t_1,\ldots,t_m$.}\end{definition}\vspace{3mm}
Let $A$ be a P.I. algebra over $F$. It is well known that $T(A)$ is a multi-homogeneous ideal of $\f$. Then, if $\overline{t}:=(t_1,\ldots,t_m),$ we denote by $$H(A,\T):=H\left({F\left<x_1,x_2,\ldots,x_m\right>}/({F\left<x_1,x_2,\ldots,x_m\right>\cap T(A)}),\T\right)$$ the Hilbert series of the relatively free algebra in $m$ variables.\vspace{3mm}

Hilbert series is related with usual operations between graded vector spaces. In fact

\begin{prop}\label{prop6}
Let V, W be $\Z_m$-graded vector spaces and U be a $\Z_m$-graded vector subspace of V. Then
\begin{itemize}
\item $H(V\oplus W,t_1,\ldots,t_m)=H(V,t_1,\ldots,t_m)+H(W,t_1,\ldots,t_m)$
\item $H(V\otimes W,t_1,\ldots,t_m)=H(V,t_1,\ldots,t_m)\cdot H(W,t_1,\ldots,t_m)$
\item $H(V,t_1,\ldots,t_m)=H(V/U,t_1,\ldots,t_m)+H(U,t_1,\ldots,t_m).$
\end{itemize}
\end{prop}

\begin{definition} Let $T$ be a tableau of shape $\lambda$ filled in with natural numbers $\{1,\ldots,k\}$ and let $d_i$ be the multiplicity of $i$ in $T$. Let $$S_{\lambda}=\sum_{\text{\rm $T_{\lambda}$ semistandard}}t^{T_{\lambda}},$$ where $t^{T_{\lambda}}=t_1^{d_1}t_2^{d_2}\cdots t_k^{d_k}.$ We say $S_{\lambda}$ is the \emph{Schur Function of $\lambda$ in the variables $t_1,\ldots,t_k.$}\end{definition}\vspace{3mm}

It is well known (see \cite{Sag}, Chapter 4) that the product of the Schur's functions $$S_{\lambda}(t_1,\ldots,t_k)S_{\mu}(t_1,\ldots,t_k)$$ corresponds in a natural way to $$(\lambda\otimes\mu)^{\uparrow S_n},$$ where $\lambda\vdash l$, $\mu\vdash m$ and $l+m=n.$ Then in the computation of products of Schur functions we are allowed to use the combinatorial tool of the Littlewwod-Richardson Rule. For a more detailed account about the Littlewood-Richardson rule, see \cite{Sag}.\vspace{3mm}

Hilbert series is strictly connected to the sequence of cocharacters of P.I. algebras. Indeed, we have the following proposition:

\begin{prop}\label{prop5}
Let A be a P.I. algebra and let $\chi_n(A)=\sum_{\lambda\vdash n}m_{\lambda}(A)\chi_{\lambda}$ its n-th cocharacter. Let $H(A,t_1,\ldots,t_k)$ the Hilbert series of A, then $$H(A,t_1,\ldots,t_k)=\sum_{n\geq0}\sum_{\lambda\in H(k,0,n)}m_{\lambda}(A)S_{\lambda}(t_1,\ldots,t_k).$$
\end{prop}\vspace{3mm}

At the light of the previous proposition, in order to compute the cocharacters of a P.I. algebra $A$, it is sufficient to argue in terms of its Hilbert series.

We can define in an analogous way the so called \emph{proper Hilbert series} of a P.I. algebra.

\begin{definition} Let $A$ be a P.I. algebra over $F$, then, if $\overline{t}:=(t_1,\ldots,t_m),$ we denote by $$H^B(A,\T):=H\left({B(x_1,\ldots,x_m)}/({B(x_1,\ldots,x_m)\cap T(A)},\T\right))$$ the Hilbert series of the relatively free algebra of proper polynomial in $m$ variables.\end{definition}\vspace{3mm}

Hilbert series and proper Hilbert series are related by the following proposition (see for example \cite{Dr03}, Theorem 4.3.12 $(i)$):

\begin{prop}\label{prop4} Let A be a P.I. algebra. Then $$H(A,\T)=\prod_{i=1}^m\frac{1}{1-t_i}\cdot H^B(A,\T).$$
\end{prop}

\newtheorem{rem}{Remark}[section]\begin{rem} The rational function $\prod_{i=1}^m\frac{1}{1-t_i}$ can be expressed in terms of Schur functions as $\sum_{k\geq0}S_{(k)}.$\end{rem}

\section{The Theorem of Lewin and Hilbert series of upper triangular matrices}

Let $I$ and $J$ be two $T$-ideals. Consider the quotient algebras
$\f/I,\f/J$ and let $U$ be a $\f/I$-$\f/J$-bimodule. We define:
$$R=\left(
    \begin{array}{cc}
      \f/I & U \\
      0 & \f/J \\
    \end{array}
  \right)
.$$
Fix ${u_i}$ a countable set of elements of $U$. Then $\varphi:x_i\rightarrow a_i$ defines an algebra homomorphism,
where:
$$a_i=\left(
       \begin{array}{cc}
         x_i+I & u_i \\
         0 & x_i+J \\
       \end{array}
     \right)
.$$
If $f(x_1,\ldots,x_n)\in\f$ one has that $f(x_1,\ldots,x_n)\rightarrow f(a_1,\ldots,a_n)$, where:
$$f(a_1,\ldots,a_n)=\left(
                      \begin{array}{cc}
                        f(x_1,\ldots,x_n)+I & \delta(f) \\
                        0 & f(x_1,\ldots,x_n)+J \\
                      \end{array}
                    \right)
$$
and $\delta(f)$ is some element of $U$. Then $IJ\subseteq ker(\varphi)=I\cap J\cap ker(\delta(f))$ and $\delta:\f\rightarrow U$ is
an $F$-derivation.\vspace{3mm}

\newtheorem{theorem}{Theorem}[section]
\begin{theorem}\label{teo1}(Lewin \cite{Lew}). If $\{u_i\}$ is a countable free set of elements of the bimodule U
then for the homomorphism $\varphi$ defined by $\{u_i\}$, we have $\ker(\varphi)=IJ.$
\end{theorem}\vspace{3mm}
\newtheorem{cor}{Corollary}[section]
\begin{cor}\label{cor}
If the bimodule U contains a countable free set $\{u_i\}$ for any i, then $T(R)=IJ.$
\end{cor}
\proof We have to prove only the "left to right" inclusion. By the Theorem of Lewin $\ker(\varphi)=IJ$ but $T(R)\subseteq\ker(\varphi)=IJ$ and we are done.\endproof

Let $A,B$ be $F$-algebras, suppose they are both P.I. and let $U$ be an $A-B$-bimodule, then we can consider $$R=\left(
       \begin{array}{cc}
         A & U \\
         0 & B \\
       \end{array}
     \right)$$
that is still an $F$-algebra and a P.I. algebra such that $T(R)\supseteq T(A)T(B)$. Suppose now that $T(R)=T(A)T(B),$ then we have a formula that relates the Hilbert series of $R$ with the Hilbert series of $A$ and $B$.

\begin{equation}\label{eq2}
H(R,\T)=H(A,\T)+H(B,\T)+\left(S_{(1)}-1\right)H(A,\T)H(B,\T).
\end{equation}

The equation (\ref{eq2}) gives us the possibility to compute the proper Hilbert series of the algebra $R$ using Proposition \ref{prop4}:

\begin{equation}\label{eq3}
H^B(R,\T)=H^B(A,\T)+H^B(B,\T)+\frac{\left(S_{(1)}-1\right)}{\prod_{i=1}^m(1-x_i)}H^B(A,\T)H^B(B,\T).
\end{equation}

\section{Cocharacters of G}

In this section, we compute the Hilbert series of $G$ that is one of the five minimal algebra of exponent 2 in the classification of Giambruno and Zaicev \cite{Giam03}. This variety has been studied by Stoyanova-Venkova in \cite{Sto}, too.

The matrix algebra
$$G=\left(
                                       \begin{array}{cc}
                                         E & E \\
                                         0 & E^0 \\
                                       \end{array}
                                     \right)
$$
is a P.I. algebra such that $T(G)=T(E)T(E^0).$ We will use the equation (\ref{eq2}) to compute its Hilbert series, then we can compute its cocharacter sequence at the light of Proposition \ref{prop5}.

Let us divide the proof in some lemmas using the well known facts that $H(E,\T)=\sum S_{(k,1^l)}$ (see, \cite{Ols}) and $H(E^0,\T)=\sum S_{(k)}.$

\newtheorem{lem}{Lemma}[section]
\begin{lem}\label{lemma9}$
H(E,\T)H(E^0,\T)=\sum m_\lambda S_\lambda,$ where $\lambda=(k_1,k_2,1^l)$ and
\begin{displaymath}\begin{array}{cc}m_\lambda=2(k_1-k_2+1) \textrm{\ if $k_2\geq1$}\\
m_\lambda=k_1+1 \textrm{\ if $k_2=l=0$}.\\
\end{array}
\end{displaymath}
\end{lem}

\proof Suppose $l\geq1$, then by the Littlewood-Richardson rule one easily has that $(k_1,k_2,1^l)$ comes from the tensor product $(k_1-i,1^{l+1})\otimes(i+k_2-1)$ or from $(k_1-i,1^{l})\otimes(i+k_2)$ for $i=0,\ldots,k_1-k_2$, so the total multiplicity is $2(k_1-k_2+1)$. We argue analogously for the other case.

\endproof\vspace{3mm}

\begin{lem}\label{lemma10}$
S_{(1)}H(E,\T)H(E^0,\T)=\sum m_\lambda S_\lambda,$ where $\lambda=(k_1,k_2,1^l)$ or $\lambda=(k_1,k_2,2,1^l).$ If $\lambda=(k_1,k_2,1^l),$ then
\begin{displaymath}\begin{array}{cc}m_\lambda=6(k_1-k_2+1) \textrm{\ if $\ k_2\geq2$ and $l\geq1$}\\
m_\lambda=4k_1-2 \textrm{\ if $k_2=0$ and $l\geq2$}\\
m_\lambda=3k_1-1 \textrm{\ if $k_2=0$ and $l=1$}\\
m_\lambda=4(k_1-k_2+1) \textrm{\ if $l=0$}\\
m_\lambda=k_1 \textrm{\ if $k_2=l=0$}.\\
\end{array}
\end{displaymath}
If $\lambda=(k_1,k_2,2,1^l),$ then \begin{displaymath}m_\lambda=2(k_1-k_2+1).\end{displaymath}
\end{lem}

\proof Suppose $\lambda=(k_1,k_2,1^l)$ and $l\geq1$, then by the Littlewood-Richardson rule one easily has that $(k_1,k_2,1^l)$ comes from the tensor product $(k_1,k_2,1^{l-1})\otimes(1),$ or $(k_1,k_2-1,1^{l})\otimes(1)$ or $(k_1-1,k_2,1^{l})\otimes(1)$. By previous Lemma, their multiplicities are respectively equal to $2(k_1-k_2+1)$, $2(k_1-k_2+2)$ and $2(k_1-k_2)$, so the total multiplicity equals $2(k_1-k_2+1)+2(k_1-k_2+2)+2(k_1-k_2)=6(k_1-k_2+1).$

Suppose now $\lambda=(k_1,k_2)$, then one has that $(k_1,k_2)$ comes from the tensor product $(k_1,k_2-1)\otimes(1),$ or $(k_1-1,k_2)\otimes(1)$. By the previous Lemma, their multiplicities are respectively equal to $2(k_1-k_2+2)$ and $2(k_1-k_2)$, so the total multiplicity equals $2(k_1-k_2+2)+2(k_1-k_2)=4(k_1-k_2+1).$

Finally, if $\lambda=(k_1,k_2,2,1^l)$ one has that $(k_1,k_2,2,1^l)$ comes from the tensor product $(k_1,k_2,1^{l+1})\otimes(1)$ only. By previous Lemma, its multiplicity equals $2(k_1-k_2+1)$ and we are done. The other cases are treated similarly.\endproof\vspace{3mm}

\begin{lem}\label{lemma11}$
(S_{(1)}-1)H(E,\T)H(E^0,\T)=\sum m_\lambda S_\lambda,$ where $\lambda=(k_1,k_2,1^l)$ or $\lambda=(k_1,k_2,2,1^l).$ If $\lambda=(k_1,k_2,1^l),$ then
\begin{displaymath}\begin{array}{cc}m_\lambda=4(k_1-k_2+1) \textrm{if $\ k_2\geq2$ and $l\geq1$}\\
m_\lambda=2k_1-2 \textrm{\ if $k_2=0$ and $l\geq2$}\\
m_\lambda=k_1-1 \textrm{\ if $k_2=0$ and $l=1$}\\
m_\lambda=2(k_1-k_2+1) \textrm{\ if $l=0$}\\
m_\lambda=-1 \textrm{\ if $k_2=l=0$}.\\
\end{array}
\end{displaymath}
If $\lambda=(k_1,k_2,2,1^l),$ then \begin{displaymath}m_\lambda=2(k_1-k_2+1).\end{displaymath}
\end{lem}

\proof We have just to use Lemmas \ref{lemma10} and \ref{lemma11}.\endproof\vspace{3mm}

Finally, the complete Hilbert series.

\begin{prop}\label{lemma12}$
H(G,\T)=\sum m_\lambda S_\lambda,$ where $\lambda=(k_1,k_2,1^l)$ or $\lambda=(k_1,k_2,2,1^l).$ If $\lambda=(k_1,k_2,1^l),$ then
\begin{displaymath}\begin{array}{cc}m_\lambda=4(k_1-k_2+1) \textrm{\ if $\ k_2\geq2$ and $l\geq1$}\\
m_\lambda=2k_1-1 \textrm{\ if $k_2=0$ and $l\geq2$}\\
m_\lambda=2(k_1-k_2+1) \textrm{\ if $l=0$}\\
m_\lambda=k_1 \textrm{\ if $l=1$}\\
m_\lambda=1 \textrm{\ if $k_2=l=0$}.\\
\end{array}
\end{displaymath}
If $\lambda=(k_1,k_2,2,1^l),$ then \begin{displaymath}m_\lambda=2(k_1-k_2+1).\end{displaymath}
\end{prop}
\proof Straightforward.\endproof

\newtheorem{ex}{Example}[section]\begin{ex} Using Proposition \ref{lemma12}, we have that
$$\chi_1(G)=(1),$$
$$\chi_2(G)=(2)+(1^2),$$
$$\chi_3(G)=(3)+2(2,1)+(1^3),$$
$$\chi_4(G)=(4)+3(3,1)+2(2^2)+3(2,1^2)+(1^4),$$
$$\chi_5(G)=(5)+4(4,1)+4(3,2)+5(3,1^2)+4(2^2,1)+3(2,1^3)+(1^5),$$
$$\chi_6(G)=(6)+5(5,1)+6(4,2)+7(4,1^2)+8(3,2,1)+2(3^2)+5(3,1^3)+2(2^3)+4(2^2,1^2)+3(2,1^4)+(1^6).$$ \end{ex}

\section{Cocharacters of $UT_2(E)$}

We compute the Hilbert series of $UT_2(E)$ starting from its proper Hilbert series and using the combinatorial properties of the Littlewood-Richardson rule.

Consider the matrix algebra $$UT_2(E)=\left(
       \begin{array}{cc}
         E & E \\
         0 & E \\
       \end{array}
     \right)$$
and let $R:=UT_2(E)$.
At the light of Theorem \ref{teo1} and its following Corollary, we have that $T(R)=T(E)T(E)$, then the Proposition (\ref{prop4}) gives us that

\begin{equation}\label{eq3}
H^B(R,\T)=2H^B(E,\T)+\frac{\left(S_{(1)}-1\right)}{\prod_{i=1}^m(1-x_i)}(H^B(E,\T))^2.
\end{equation}

It is convenient to break the computation of $H^B(R,\T)$ into some lemmas.

\begin{lem}\label{lemma1}
$(H^B(E,\T))^2=\sum m_\lambda S_\lambda,$ where $$\lambda=\left(2^{\mu_2},1^{\mu_1}\right)$$ and \begin{displaymath}m_\lambda=\left\{ \begin{array}{ll}
\frac{\mu_1-\mu_2}{2}+1 & \textrm{if $\mu_2$ even} \\
\frac{\mu_1-\mu_2}{2} & \textrm{if $\mu_2$ odd}.
\end{array} \right.
\end{displaymath}
\end{lem}

\proof It is well known that $H^B(E,\T)=\sum_{k\geq0} S_{(1^{2k})}$ (see \cite{Dr03}, Chapters 4, 12), so let us check out the multiplicities of $\left[(\sum_{k\geq0}S_{(2k)})^2\right]'.$ By Littlewood-Richardson rule we have that the only allowed partitions in the tensor product $$\yng(6)\otimes\yng(4)$$ are of the type $\mu=(\mu_1,\mu_2)$ as in the picture below $$\mu=\yng(7,3).$$ Due to the parity of $|\mu|,$ we have that $\mu_1,\mu_2$ are both even or both odd numbers. If $\mu_2$ is even we have exactly $\frac{\mu_1-\mu_2}{2}+1$ allowed partitions. If $\mu_2$ is odd, $\mu$ doesn't occur as $(\mu_1)\otimes(\mu_2)$ so we have exactly $\frac{\mu_1-\mu_2}{2}$ allowed partitions.\endproof\vspace{3mm}

\begin{rem} The partitions $\lambda_1=(2,1^l),$ if $l$ is odd and $\lambda_2=(2)$ are not allowed in the previous decomposition. In fact, $|\lambda_1|$ is odd and $\lambda_2$ comes only from $1\otimes(2)$ or $(1)\otimes(1)$ and both of $(1)$ and $(2)$ do not appear in the decomposition of $H^B(E,\T)$.\end{rem}

\begin{lem}\label{lemma2}
$\frac{1}{\prod_{i=1}^m(1-x_i)}(H^B(E,\T))^2=\sum m_\lambda S_{\lambda},$ where
\begin{displaymath}
\begin{array}{c}{\lambda}=(k,2^m,1^l)\\
\end{array}
\end{displaymath}
and
\begin{displaymath}
\begin{array}{cc}
m_{{\lambda}}=\left\{\begin{array}{cc} \frac{l}{2}+1 & \textrm{if l is even} \\
                                                 \frac{l-1}{2}+1 & \textrm{if l is odd}\end{array} \right. if $\ k=m=0.$\\
m_{{\lambda}}=l+1 \textrm{\ otherwise.}\\\vspace{3mm}
\end{array}
\end{displaymath}
\end{lem}

\proof We argue only for $\lambda=(k,2^m,1^l),$ with $m\geq1$, because the case $m=0$ is treated similarly. By the Littlewood-Richardson rule one has that
$(k,2^m,1^l)$ occurs in $(2^{m+1},1^{l})\otimes(k-2)$, $(2^m,1^{l})\otimes(k)$ if $l$ is even, $(2^{m+1},1^{l-1})\otimes(k-1),$ $(2^{m},1^{l+1})\otimes(k-1)$ if $l$ is odd, then by Lemma \ref{lemma1}, $(2^{m+1},1^{l})$, $(2^{m},1^{l})$, $(2^{m+1},1^{l-1})$ and $(2^{m},1^{l+1})\otimes(k-1),$ have multiplicities respectively equal to $\frac{l}{2},\frac{l}{2}+1,\frac{l-1}{2}$, and $\frac{l+1}{2}+1$ if $m$ is even, equal to $\frac{l}{2}+1,\frac{l}{2},\frac{l-1}{2}+1$, and $\frac{l+1}{2}$ if $m$ is odd. Then $m_{\lambda}=\frac{l}{2}+\frac{l}{2}+1=l+1$ if $l$ is even and $m$ is even;
$m_{\lambda}=\frac{l}{2}+1+\frac{l}{2}=l+1$ if $l$ is even and $m$ is odd;
$m_{\lambda}=\frac{l-1}{2}+\frac{l+1}{2}+1=l+1$ if $l$ is odd and $m$ is even;
$m_{\lambda}=\frac{l-1}{2}+1+\frac{l+1}{2}=l+1$ if $l$ is odd and $m$ is odd, from which the assertion.

If $\lambda=(1^l),$ By the Littlewood-Richardson rule one has that
$(1^l)$ occurs in $(1^{l})\otimes1$ if $l$ is even, $(1^{l-1})\otimes(1)$ if $l$ is odd, then by Lemma \ref{lemma1}, it has multiplicities respectively equal to $\frac{l}{2}+1$ and $\frac{l-1}{2}+1$ and we are done.\endproof\vspace{3mm}

\begin{lem}\label{lemma3}
$\frac{S_{(1)}}{\prod_{i=1}^m(1-x_i)}(H^B(E,\T))^2=\sum m_\lambda S_{\lambda},$ where
\begin{displaymath}
\begin{array}{c}{\lambda}=(k,2^m,1^l)\\
\textrm{or }{\lambda}=(k,3,2^m,1^l).\\
\end{array}
\end{displaymath}
If $\lambda=(k,2^m,1^l),$ then
\begin{displaymath}
\begin{array}{cc}m_{{\lambda}}=3(l+1) \textrm{ if $\ k\geq3,m\geq1$}\\
m_{{\lambda}}=2(l+1) \textrm{ if $k=2,\ m\geq1$}\\
m_{{\lambda}}=2l+1 \textrm{ if $\ k\geq3,m=0$}\\
m_{{\lambda}}=\left\{\begin{array}{cc} \frac{l}{2} & \textrm{if l is even} \\
                                                 \frac{l-1}{2}+1 & \textrm{if l is odd}\\
                               \end{array} \right.\textrm{ if $m=k=0$}\\
m_{{\lambda}}=\left\{\begin{array}{cc} l+\frac{l}{2}+1 & \textrm{if l is even} \\
                                                 l+\frac{l+1}{2}+1 & \textrm{if l is odd}\\
                               \end{array} \right.\textrm{ if $k=2,m=0$}\\
\end{array}
\end{displaymath}
If $\lambda=(k,3,2^m,1^l),$ then \begin{displaymath}m_{{\lambda}}=l+1.\end{displaymath}
\end{lem}

\proof If $\lambda=(k,2^m,1^l)$ and $k\geq3,m\geq1,$ by the Littlewood-Richardson rule one has that
$\lambda$ occurs in $(k,2^{m},1^{l-1})\otimes(1)$, $(k,2^{m-1},1^{l+1})\otimes(1)$, $(k-1,2^{m},1^{l})\otimes(1)$ and by Lemma \ref{lemma2}, their multiplicities are $l$, $l+2$, $l+1$ respectively. Then the total multiplicity of $(k,2^m,1^l)$ is $l+l+2+l+1=3l+3=3(l+1)$.

If $\lambda=(2^m,1^l)$ and $m\geq1,$ by the Littlewood-Richardson rule one has that
$\lambda$ occurs in $(2^{m},1^{l-1})\otimes(1)$, $(2^{m-1},1^{l+1})\otimes(1)$, and by Lemma \ref{lemma2}, their multiplicities are $l$, and $l+2$ respectively. Then the total multiplicity of $(2^m,1^l)$ is $l+l+2=2(l+1).$

If $\lambda=(k,1^l)$ one has that
$\lambda$ occurs in $(k,1^{l-1})\otimes(1)$, $(k-1,1^{l})\otimes(1)$, and by Lemma \ref{lemma2}, their multiplicities are $l$, and $l+1$ respectively. Then the total multiplicity of $(2^m,1^l)$ is $l+l+1=2l+1.$

If $k=m=0$, one has that if $l$ is even, $l$ comes from $(1^{l-1})\otimes(1)$, then $l-1$ is odd and by previous Lemma, its multiplicity is $\frac{l-1-1}{2}+1=\frac{l}{2};$ if $l$ is odd, $l$ comes from $(1^{l-1})\otimes(1)$, then $l-1$ is even and by previous Lemma, its multiplicity is $\frac{l-1}{2}+1.$

Finally, if If $\lambda=(k,3,2^m,1^l)$, one has that
$\lambda$ occurs in $(k,2^{m+1},1^{l})\otimes(1)$ only and its multiplicity is $l+1$. The other cases are treated similarly.\endproof\vspace{3mm}

\begin{lem}\label{lemma4}
$\frac{S_{(1)}-1}{\prod_{i=1}^m(1-x_i)}(H^B(E,\T))^2=\sum m_\lambda S_{\lambda},$ where
\begin{displaymath}
\begin{array}{c}{\lambda}=(k,2^m,1^l)\\
\textrm{or }{\lambda}=(k,3,2^m,1^l).\\
\end{array}
\end{displaymath}
If $\lambda=(k,2^m,1^l),$ then
\begin{displaymath}
\begin{array}{cc}m_{{\lambda}}=2(l+1) \textrm{ if $\ k\geq3,m\geq1$}\\
m_{{\lambda}}=l+1 \textrm{ if $k=2,\ m\geq2$}\\
m_{{\lambda}}=l \textrm{ if $\ k\geq3,m=0$}\\
m_{{\lambda}}=\left\{\begin{array}{cc} -1 & \textrm{if l is even} \\
                                                 0 & \textrm{if l is odd}\\
                               \end{array} \right.\textrm{ if $m=k=0$}\\
m_{{\lambda}}=\left\{\begin{array}{cc} \frac{l}{2} & \textrm{if l is even} \\
                                                 \frac{l+1}{2} & \textrm{if l is odd}\\
                               \end{array} \right.\textrm{ if $k=2,m=0$}\\
\end{array}
\end{displaymath}
If $\lambda=(k,3,2^m,1^l),$ then \begin{displaymath}m_{{\lambda}}=l+1.\end{displaymath}
\end{lem}

\proof We have just to combine Lemmas \ref{lemma2} and \ref{lemma3}.\endproof\vspace{3mm}

Now we have all the results to state the following proposition:
\begin{prop}\label{lemma5}
$H^B(R,\T)=2H^B(E,\T)+\frac{S_{(1)}-1}{\prod_{i=1}^m(1-x_i)}(H^B(E,\T))^2=\sum m_\lambda S_{\lambda},$ where
\begin{displaymath}
\begin{array}{c}{\lambda}=(k,2^m,1^l)\\
\textrm{or }{\lambda}=(k,3,2^m,1^l).\\
\end{array}
\end{displaymath}
If $\lambda=(k,2^m,1^l),$ then
\begin{displaymath}
\begin{array}{cc}m_{{\lambda}}=2(l+1) \textrm{ if $\ k\geq3,m\geq1$}\\
m_{{\lambda}}=l+1 \textrm{ if $\ m\geq2$}\\
m_{{\lambda}}=l \textrm{ if $\ k\geq3,m=0$}\\
m_{{\lambda}}=\left\{\begin{array}{cc} 1 & \textrm{if l is even} \\
                                                 0 & \textrm{if l is odd}\\
                               \end{array} \right.\textrm{ if $m=k=0$}\\
m_{{\lambda}}=\left\{\begin{array}{cc} \frac{l}{2} & \textrm{if l is even} \\
                                                 \frac{l+1}{2} & \textrm{if l is odd}\\
                               \end{array} \right.\textrm{ if $k=2,m=0$}\\
\end{array}
\end{displaymath}
If $\lambda=(k,3,2^m,1^l),$ then \begin{displaymath}m_{{\lambda}}=l+1.\end{displaymath}
\end{prop}

\proof We have to add 2 at the multiplicity of $(1^l)$ where $l$ is even, computed in the previous lemma.\endproof\vspace{3mm}

Now, we have all informations about proper Hilbert series of the algebra $R$ but, at the light of Proposition \ref{prop4}, we have that $$H(R,\T)=\frac{1}{\prod_{i=1}^m1-t_i}H^B(R,\T),$$ so we can state the following proposition:
\begin{prop}\label{prop1}
$H(R,\T)=\sum m_\lambda S_{\lambda},$ where
\begin{displaymath}
\begin{array}{c}{\lambda}=(k_1,k_2,2^m,1^l)\textrm{ or ${\lambda}=(k_1,k_2,3,2^m,1^l)$}.\\
\end{array}
\end{displaymath}
If ${\lambda}=(k_1,k_2,2^m,1^l),$ then
\begin{displaymath}
\begin{array}{cc}m_{{\lambda}}=12(k_1-k_2+1)(l+1) \textrm{\ \ \ if $k_1\geq k_2\geq3,\ m\geq1$}\\
m_{{\lambda}}=4(k_1-k_2+1)(2l+1) \textrm{\ \ \ if $k_1\geq k_2\geq3,\ m=0$}\\
m_{{\lambda}}=8(k_1-2)(l+1)+4(l+1) \textrm{\ \ \ if $k_1\geq k_2=2,\ m\geq1$}\\
m_{{\lambda}}=3(k_1-2)(2l+1)+3l+2 \textrm{\ \ \ if $k_1\geq k_2=2,\ m=0$}\\
m_{{\lambda}}=(k_1-2)(2l-1)+l+1 \textrm{\ \ \ if $k_1\geq2,\ k_2=0,\ m=0,\ l\geq1$}\\
m_{{\lambda}}=1 \textrm{\ \ \ if $\lambda=(1^l)$ or $\lambda=(k)$}.\\
\end{array}
\end{displaymath}
If ${\lambda}=(k_1,k_2,2^m,1^l),$ then
\begin{displaymath}
\begin{array}{cc}m_{{\lambda}}=4(k_1-k_2+1)(l+1) \textrm{\ if $k_2\geq3m\geq1$}.\\
\end{array}
\end{displaymath}

\end{prop}

\proof Suppose firstly $k_2\geq3,m\geq1$, then by Proposition \ref{prop4} and the following Remark, we have that $$m_{{\lambda}}=\sum_{i=0}^{k_1-k_2}m_{(k_1-i,3,2^m,1^l)}+m_{(k_1-i,3,2^{m-1},1^{l+1})}+
m_{(k_1-i,3,2^m,1^{l-1})}+m_{(k_1-i,3,2^{m-1},1^{l})}+$$ $$+\sum_{i=0}^{k_1-k_2}
m_{(k_1-i,2^{m+1},1^l)}+m_{(k_1-i,2^{m},1^{l+1})}+m_{(k_1-i,2^{m+1},1^{l-1})}+m_{(k_1-i,2^{m},1^l)}$$ and, by Proposition \ref{lemma5}, $m_{\lambda}=(k_1-k_2+1)(l+1+l+2+l+l+1)+(k_1-k_2+1)(2(l+1)+2(l+2)+2l+2(l+1))=4(k_1-k_2+1)(l+1)+8(k_1-k_2+1)=12(k_1-k_2+1)(l+1).$

Let $\lambda=(k)$, then we have that $$m_{{\lambda}}=\sum_{i=0}^{k-1}
m_{(k-i)}+m_{1}$$ and, by Proposition \ref{lemma5},
$m_{\lambda}=1$.

Finally, Let $\lambda=(1^l)$, then we have that $$m_{{\lambda}}=
m_{(1^l)}+m_{(1^{l-1})}$$ and, by Proposition \ref{lemma5},
$m_{\lambda}=1$. The other cases are treated similarly.\endproof
\vspace{3mm}

\begin{ex}Using Proposition \ref{prop5} and the previous proposition, we have that
$$\chi_1(R)=(1),$$
$$\chi_2(R)=(2)+(1^2),$$
$$\chi_3(R)=(3)+2(2,1)+(1^3),$$
$$\chi_4(R)=(4)+3(3,1)+2(2^2)+3(2,1^2)+(1^4),$$
$$\chi_5(R)=(5)+4(4,1)+5(3,2)+6(3,1^2)+5(2^2,1)+4(2,1^3)+(1^5),$$
$$\chi_6(R)=(6)+5(5,1)+8(4,2)+9(4,1^2)+14(3,2,1)+4(3^2)+9(3,1^3)+(2^3)+8(2^2,1^2)+5(2,1^4)+(1^6).$$ \end{ex}

\section{$\Z_2$-graded cocharacters of $UT_2(E)$}

 We compute directly the $\Z_2$-graded cocharacters of $UT_2(E)$ starting from the cocharacters of $UT_2(F)$, then we will use a formula due to Di Vincenzo and Nardozza (see \cite{Divi03}).

Using direct computations, is easy to prove (indeed, see \cite{Ber02}) that the ordinary cocharacters sequence of $$S=\left(
       \begin{array}{cc}
         F & F \\
         0 & F \\
       \end{array}
     \right)$$
is the following

\begin{prop}\label{prop2}
For any $n\geq0$, $\chi_n(S)=\sum_{i=1}^{3}m_\lambda^{(i)} S_{\lambda^{(i)}},$ where
\begin{displaymath}
\begin{array}{c}{\lambda}^{(1)}=(n)\\
{\lambda}^{(2)}=(k_1,k_2)\\
{\lambda}^{(3)}=(k_1,k_2,1)\\
\end{array}
\end{displaymath}
and
\begin{displaymath}
\begin{array}{cc}m_{{\lambda}^{(1)}}=1\\
m_{{\lambda}^{(2)}}=k_1-k_2+1\\
m_{{\lambda}^{(3)}}=k_1-k_2+1.\\
\end{array}
\end{displaymath}
\end{prop}\vspace{3mm}

The algebra $R$ has a natural structure of $\Z_2$-graded algebra, where $$R^0=\left(
       \begin{array}{cc}
         E^0 & E^0 \\
         0 & E^0 \\
       \end{array}
     \right)$$
 and $$R^1=\left(
       \begin{array}{cc}
         E^1 & E^1 \\
         0 & E^1 \\
       \end{array}
     \right),$$
and $E^0,E^1$ are, respectively, the 0 and the 1 part of the natural $\Z_2$-grading of $E$. As a $\Z_2$-graded algebra, $R$ is naturally isomorphic to $S\otimes E$. We have the following proposition due to Di Vincenzo and Nardozza:

\begin{prop}\label{prop3}
Let $k,l\in\N$ such that $k+l=n$ and consider $H=S_k\times S_l$. If $(\chi_n(S))_{\downarrow H}=\sum m_{\lambda,\mu}\lambda\otimes\mu$, then $$\chi_n^{\Z_2}(R)=\sum_{k+l=n}\sum_{\begin{array}{c}
                                      \lambda\vdash k \\
                                      \mu\vdash l
                                    \end{array}
} m_{\lambda,\mu}\lambda\otimes\mu'.$$
\end{prop}\vspace{3mm}

At the light of Proposition \ref{prop3}, it suffices to know $(\chi_n(S))_{\downarrow H}=\sum m_{\lambda,\mu}\lambda\otimes\mu$, then we have immediately the $n$-th $\Z_2$-graded cocharacter of $UT_2(E)$. At this purpose we have the following Proposition.

\begin{prop}\label{prop7}
For any $n\geq0$, $(\chi_n(S))_{\downarrow H}=\sum m_{\lambda,\mu}\lambda\otimes\mu,$ where
\begin{displaymath}
\begin{array}{c}{\lambda}=(\lambda_1,\lambda_2,\lambda_3), \quad  \textrm{$\lambda_3\leq1$}\\
{\mu}=(\mu_1,\mu_2,\mu_3), \quad  \textrm{$\mu_3\leq1$}.\\
\end{array}
\end{displaymath}
More precisely,
$$(n)_{\downarrow H}=(k)\otimes(l)$$
$$(k_1,k_2)_{\downarrow H}=\left\{\begin{array}{cc}
                                                                  (k)\otimes(l)  \\
                                                                  (\lambda_1,\lambda_2)\otimes(l) \\
                                                                  (\lambda_1,\lambda_2)\otimes(\mu_1,\mu_2)
                                                                \end{array}\right.$$

$$(k_1,k_2,1)_{\downarrow H}=\left\{\begin{array}{cc}
                                                                  (\lambda_1,\lambda_2)\otimes(l) \\
                                                                  \begin{array}{cc}2((\lambda_1,\lambda_2)\otimes(\mu_1,\mu_2)) & \textrm{if $\mu_1-1\geq\mu_2$}\\
                                                                  (\lambda_1,\lambda_2)\otimes(\mu_1,\mu_2) & \textrm{if $\mu_1-1<\mu_2$}\\
                                                                  \end{array}\\
                                                                  (\lambda_1,\lambda_2,1)\otimes(l) \\
                                                                  (\lambda_1,\lambda_2,1)\otimes(\mu_1,\mu_2)
                                                                \end{array}\right.$$

\end{prop}

\proof The result follows from a straightforward computation using the decomposition of $\chi_n(S)$ given in Proposition \ref{prop2}. More precisely, it follows by Branching Rule (see \cite{Sag}, Chapter 2) that when we restrict the representation $\nu=\sum(k_1,k_2,k_3)$ of $S_n$, with $k_3\leq1$, to its subgroup $H$, then its $H$-irreducible components are $\lambda\otimes\mu$, where $\lambda=(k_1',k_2',k_3')$ and $\mu=(l_1,l_2,l_3)$ are such that $\nu$ appears in the tensor product $(\lambda\otimes\mu)^{\uparrow S_n}$. By Frobenius Multiplicity Law (see \cite{Sag}), the multiplicity of $\lambda\otimes\mu$ in the previous decomposition equals the multiplicity of $\nu$ in the induced representation $(\lambda\otimes\mu)^{\uparrow S_n}.$
We will argue for the irreducible cocharacters of $\chi_n(F)$, i.e. $(n)$, $(k_1,k_2)$ where $k_2\geq1$ and $(k_1,k_2,1)$. In the first case, it is easy to see that $(n)_{\downarrow H}=\m((l)\otimes(k)),$ where $\lambda=(k)$ and $\mu=(l)$. By the Littlewood-Richardson Rule, the multiplicity of $(n)$ in the induced representation $(\lambda\otimes\mu)^{\uparrow S_n}$ is 1.
Let $\nu=(k_1,k_2)$. Then $$(k_1,k_2)_{\downarrow H}=\left\{\begin{array}{cc}
                                                                  \m((k)\otimes(l)) & a) \\
                                                                  \m((\lambda_1,\lambda_2)\otimes(l)) & b) \\
                                                                  \m((\lambda_1,\lambda_2)\otimes(\mu_1,\mu_2)) & c) \\
                                                                \end{array}\right.$$
The case $a)$ has been yet treated. Consider the case $b)$. Here, the only way to obtain $\nu$ is adding $k_2-\lambda_2$ $\young(1)$ to $\lambda_2$ so the multiplicity of $\nu$ in $(\lambda_1,\lambda_2)\otimes(l)$ is 1. Even in the case $c)$, the only possible way to obtain $\nu$ is adding all the boxes $\young(2)$ to $\lambda_2$ so the multiplicity of $\nu$ in $(\lambda_1,\lambda_2)\otimes(l)$ is still 1.
Finally, let $\nu=(k_1,k_2,1)$. Then $$(k_1,k_2,1)_{\downarrow H}=\left\{\begin{array}{cc}
                                                                  \m((\lambda_1,\lambda_2)\otimes(l)) & a) \\
                                                                  \m((\lambda_1,\lambda_2)\otimes(\mu_1,\mu_2)) & b) \\
                                                                  \m((\lambda_1,\lambda_2,1)\otimes(l)) & c) \\
                                                                  \m((\lambda_1,\lambda_2,1)\otimes(\mu_1,\mu_2)) & d)
                                                                \end{array}\right.$$
Firstly, we note that the cases $a),c),d)$ are similar to those computed for $(k_1,k_2)_{\downarrow H}$. Thus we have to argue only for the case $b)$. Suppose $\mu_1-1\geq\mu_2$, then we have to add a final box $\young(1)$ or $\young(2)$ to $\lambda$ if $\mu_1-1\geq\mu_2$, finally we have to add all the remaining $\young(2)$ in the only possible way. If $\mu_1-1<\mu_2$, we have to add only the final box $\young(2)$ to $\lambda$, finally we have to add all the remaining $\young(2)$ in the only possible way, from which the assertion.

\endproof

\section{Acknowledgements}

The author is very grateful to his advisor, Prof. Onofrio Mario Di Vincenzo and to Prof. Vesselyn Drensky for the useful discussions and suggestions.

\end{document}